\documentclass[12pt]{article}
\usepackage{amsmath,amssymb,amsthm}
\usepackage{fullpage}

\newtheorem{rmk}{Remark}
\newtheorem{thm}{Theorem}
\newtheorem{cor}{Corollary}
\newtheorem{deff}{Definition}
\newcommand{\dis}{\displaystyle}
\numberwithin{equation}{subsection}

\begin{document}

\title{\textsc{Descartes}' Perfect Lens}
\author{Mark B. Villarino\\
Depto.\ de Matem\'atica, Universidad de Costa Rica,\\
2060 San Jos\'e, Costa Rica}
\date{\today}

\maketitle

 \begin{abstract}
We give a new, elementary,  purely analytical development of \textsc{Descartes}' theorem that  a smooth connected surface is a perfect focusing lens if and only if it is a connected subset of the ovoid obtained by revolving a cartesian oval around its axis of symmetry.\end{abstract}

\tableofcontents

\section{Introduction}

Almost two thousand three hundred years ago, the hellenistic mathematician \textsc{Diokles} \cite{Dio} gave the first proof that a mirror in the shape of a paraboloid of revolution reflects all incident light rays, which are parallel to its axis of symmetry, to a single point, which \textsc{Kepler} \cite{Kep}, in 1604,  called the \emph{focus}.

After the advent of the calculus it was possible to prove that the \emph{only} such reflecting surface is generated by revolving a (proper or degenerate) parabola around its axis of symmetry.  This is a very famous and well-known result, and is treated in many easily accessible sources.  See, for example,  \textsc{Spiegel} \cite{Spi}. 

All these proofs are based on \textsc{Heron}'s \emph{\textbf{Law of Reflection}}, $\theta_I=\theta_{Rf}$, where $\theta_{Rf}$ is the angle the \emph{reflected} ray makes with the normal to the reflecting surface at the point of incidence of the incoming ray and  $\theta_I$ is the angle the\emph{ incident} ray makes with the normal.  One transforms \textsc{Heron}'s equation into the \textbf{\emph{ordinary differential equation (ODE)}} of the cross-section curve of the mirror.   \textsc{Drucker}'s paper \cite{Dru} would seem to be the final word on the subject.

Unfortunately and surprisingly, the corresponding result for a \emph{\textbf{lens}}, instead of a mirror, is  less well-known, at least among mathematicians (although \cite{Mae} is pleasant attempt to alter that). Yet the case of a lens, too, is quite fascinating and is treatable by elementary means.  The purpose of this paper is to remedy the situation and fill this gap.  

Indeed, it all started in 1637, when \textsc{Descartes} \cite{Des} asked for the \emph{refractive} analogue of the parabolic mirror:\begin{quote} \emph{Which shape of \textbf{lens} will focus all rays from one radiant point source to one single image point?} \end{quote}

We will call such a lens a \textbf{\emph{perfect lens}}.

\textsc{Descartes} discovered that the cross-section curve of the perfect lens, assumed to be a surface of revolution,  is a fourth degree curve known today as the \textbf{\emph{cartesian oval}}.  It can be defined as the locus of points the "weighted" sum of whose distances from two fixed points is a constant:
\begin{equation}
\label{car1}
\fbox{$\dis d_1+nd_2=c$}
\end{equation}where $d_1$ and $d_2$ are the distances from any point on the curve to the two fixed points, called the \emph{foci}, and $n$ is a constant.   If one focus is at the origin and the other is at the point $(b,0)$ where $b\geqslant 0$ , the equation can be written:
\begin{equation}
\label{car2}
\fbox{$\dis \left[(1-n^2)(x^2+y^2)+2n^2bx+c^2-n^2b^2\right]^2=4c^2(x^2+y^2) $}
\end{equation}

If $n=\pm 1$, the oval is the \textbf{\emph{conic section}}:
\begin{equation}
\label{car3}
\fbox{$\dis \frac{\left(x-\frac{b}{2}\right)^2}{\left(\frac{c}{2}\right)^2}-\frac{y^2}{\left(\frac{c^2-b^2}{4}\right)}=1 $}
\end{equation}

More information on cartesian ovals can be found in \cite{Ovl} and \cite{Sal} and \cite{Wei}.

\textsc{Descartes}' own treatment, which is not altogether easy to read (see \cite{BG}), shows that the oval is a solution, but does not show that it is the \emph{only} solution.

The only treatments of \textsc{Descartes}' result that we have seen in the literature do \emph{not} appear in books on mathematics (!), but rather on \emph{optics} (see \textsc{Hecht} \cite{Hec} and \textsc{Klein} \cite{Kle}) and use \textsc{Fermat}'s \emph{Principle:} \emph{A light ray traverses the path between two points which takes the \textbf{least time}.}

A non-trivial computation, based on the calculus of variations,  shows that the time the light ray takes to go from the radiant point to the image point is \emph{\textbf{constant}} for every point of the cross-section curve of the perfect lens (since if the time were different in two points of the curve, it would not be minimal), and therefore its equation is that of the cartesian oval.  

Moreover such treatments make physical assumptions about the velocity of light in different media, while, as our treatment will show, the problem is really one in pure mathematics.

We have \emph{not} seen any treatment of the subject which is founded  purely on \textbf{\emph{\textsc{Snell}'s law of refraction},}  which describes the relationship between the angle of incidence and the angle of refraction when light passes the boundary between two isotropic media (media in which the path of a light ray is a straight line).  The law states: \begin{quote}

 \emph{If $\theta_I$ is the angle the incident ray makes with the normal to the boundary at the point of refraction, and if $\theta_R$ is the angle the refracted ray makes with the normal, then at all points of the boundary the ratio\begin{equation*}
\fbox{$\dis \frac{\sin \theta_I}{\sin \theta_R}=n$}
\end{equation*}where $n$, called \emph{the index of refraction}, is \emph{\textbf{constant}} .}\end{quote}

Such a treatment of \textsc{Descartes}' theorem would seem desireable, since it is the immediate generalization of the corresponding treatment of the \textbf{\emph{perfect reflective mirror}}.

In this paper, we will present a new, self-contained,  elementary,  purely analytical proof, based on \textsc{Snell}'s Law and \textsc{Drucker}'s paper \cite{Dru},  of the following complete form of \textsc{Descartes}' theorem:


\begin{thm}(\textbf{\textsc{Descartes}' Theorem}) A smooth connected surface is a \textbf{perfect lens} if and only if it is a connected subset of the ovoid obtained by revolving a cartesian oval around its axis of symmetry. \end{thm}

\hfill$\Box$

\section{The Analytical Problem (in two dimensions)}

We begin by  solving the following two-dimensional purely analytical problem:

\begin{quote} \emph{It is required to find the equation, $f(x,y)=0$, of a smooth connected curve, $\mathbf{C}$, for which the straight lines from two fixed points cut the normal in two angles whose sines are in \textbf{constant ratio}.}\end{quote}Please note the absence of physical modeling.  The problem is purely mathematical, as its its solution.

We will find and solve an \textbf{\emph{ordinary differential equation}} (ODE) for which the equation of the curve is the general solution.  The ODE, in fact, will be a restatement of \textsc{Snell}'s law.

\begin{deff}

We call any curve $\mathbf{C}$ that solves the problem a \emph{\textbf{perfect two-dimensional lens}} with respect to the points $F$ and $F'$.
\end{deff} 

\subsection{Both Fixed Points Are Finite}
\subsubsection{The Differential Equation}
We assume a cartesian coordinate system in the xy-plane.

Let the two fixed points be O$(0,0)$ and B$(b,0)$ with $b>0$ (Here is we use the assumption that $F$ and $F'$ are finite and distinct). Let P$(x,y)$ be a variable point on the curve $f(x,y)=0$. We assume P is in the first quadrant and we assume that the curve is concave downwards at P (that is, if $y(x)$ is the function defined implicitly by the equation $f(x,y)=0$, then $y''(x_0)<0$) .  Let $l_1$ be the length of the line segment $ \overline{OP}$ and let $l_2$ be the length of  $\overline{PB}$.  Let $\overline{MN}$ be the normal to the curve where N is on the concave side of the curve and P is between M and N. Let $\theta_1:=M\widehat{P}O$, the angle that $\overline{OP}$ forms with the normal $ \overline{MN}$, and let $\theta_2:=B\widehat{P}N$, the angle that $\overline{BP}$ forms with the normal $\overline{MN}$.  Let $\overline{PT}$ be the tangent line to the curve at P where T is the point on the x-axis where the tangent line crosses it.  Let $\phi:=P\widehat{T}B$, the angle, measured counter clockwise, the tangent line forms with the x-axis.  

We will use the geometry of the figure to obtain formulas for $\sin \theta_1$ and $\sin \theta_2$ in terms of $x$, $y$, and the derivative, $y'$.  When we substitute these expressions into \textsc{Snell}'s Law, we obtain the desired\emph{ differential equation for the curve }$f(x,y)=0$.

By the law of cosines applied twice to the $\triangle OPB$
\begin{align*}
\cos P\widehat{O}B&=\dfrac{b^2+l_{1}^2-l_{2}^2}{2l_1b} & 
      \cos P\widehat{B}O&:=\dfrac{b^2+l_{2}^2-l_{1}^2}{2l_2b} \\
\end{align*}But
\begin{align*}
  \cos P\widehat{O}B&=\cos(\theta_1+\phi-90)   &
      \cos P\widehat{B}O&=\cos(90-\theta_2-\phi)\\
    &  =\sin(\theta_1+\phi)&
    &=\sin(\theta_2+\phi)\\
\end{align*}So, we obtain our fundamental formulas:
\begin{align}
\label{sin}
\sin(\theta_1+\phi)&=\dfrac{b^2+l_{1}^2-l_{2}^2}{2l_1b} & 
     \sin(\theta_2+\phi)&:=\dfrac{b^2+l_{2}^2-l_{1}^2}{2l_2b} 
\end{align}Moreover, it is evident that
\begin{align}
\label{angle}
\frac{\pi}{2} &<\theta_1+\phi<\pi&   
      0&\leqslant \theta_2+\phi<\frac{\pi}{2}
\end{align}
Now
\begin{align*}
 \dfrac{b^2+l_{1}^2-l_{2}^2}{2l_1b}&=\frac{2bx}{2l_1b}   &
      \dfrac{b^2+l_{2}^2-l_{1}^2}{2l_2b}&=\frac{2b^2-2bx}{2l_2b}\\
    &  =\frac{x}{l_1}
    &
    &=\frac{b-x}{l_2}\\
\end{align*}Therefore, equations \eqref{sin} become
\begin{align*}
\sin(\theta_1+\phi)&=  \frac{x}{l_1}&
      \sin(\theta_2+\phi)=\frac{b-x}{l_2} \\
\end{align*}By \eqref{sin} and \eqref{angle}  and the definition of the $\arcsin$ function,  we obtain
\begin{align*}
\pi-(\theta_1+\phi)&=  \arcsin\left(\frac{x}{l_1}\right)& 
      \theta_2+\phi=\arcsin\left(\frac{b-x}{l_2}\right) \\
\end{align*}and therefore
\begin{align*}
\theta_1&=  (\pi-\phi)-\arcsin\left(\frac{x}{l_1}\right)&
      \theta_2=\arcsin\left(\frac{b-x}{l_2}\right)-\phi \\
\end{align*}whence,
\begin{align*}
\sin\theta_1&=\sin(\pi-\phi)\cos\left\{\arcsin\left(\frac{x}{l_1}\right)\right\}-\cos(\pi-\phi)\sin\left\{\arcsin\left(\frac{x}{l_1}\right)\right\} \\
&=\sin \phi\sqrt{1-\frac{x^2}{l_{1}^2}}+\cos \phi\left(\frac{x}{l_1}\right)\\
&=\frac{y'}{\sqrt{1+y'^2}}\sqrt{1-\frac{x^2}{l_{1}^2}}+\frac{1}{\sqrt{1+y'^2}}\left(\frac{x}{l_1}\right)
       \end{align*}and
\begin{align*}
\sin\theta_2&=\cos \phi\sin\left \{\arcsin\left(\frac{b-x}{l_2}\right)\right\}-\sin \phi\cos \left\{\arcsin\left(\frac{b-x}{l_2}\right)\right\} \\
&=\frac{1}{\sqrt{1+y'^2}}\left(\frac{b-x}{l_2}\right)-\frac{y'}{\sqrt{1+y'^2}}\sqrt{1-\left(\frac{b-x}{l_2}\right)^2}\\
&=\frac{1}{\sqrt{1+y'^2}}\left(\frac{b-x}{l_2}\right)-\frac{y'y}{l_2\sqrt{(1+y'^2)}}.       \end{align*}   

Now, our assumption is that  \textsc{Snell}'s Law holds, i.e., that$$ \frac{\sin \theta_1}{\sin \theta_2}=n $$where $n$ is a constant,\emph{ holds for every point $P(x,y)$ of the curve.}   Substituting our two formulas for $\sin \theta_1$ and $\sin \theta_2$ into this equation gives us the equation:
\begin{align} 
\label{Snell1}
\dfrac{\dfrac{y'}{\sqrt{1+y'^2}}\sqrt{1-\dfrac{x^2}{l_{1}^2}}+\dfrac{1}{\sqrt{1+y'^2}}\left(\dfrac{x}{l_1}\right)}{\dfrac{1}{\sqrt{1+y'^2}}\left(\dfrac{b-x}{l_2}\right)-\dfrac{y'y}{l_2\sqrt{(1+y'^2)}}}=n\end{align} Solving this equation \eqref{Snell1} for $y'$ we obtain \textbf{\emph{the differential equation of the curve}}:  
\begin{equation}
\label{ODE1}
\fbox{$\dis y'=\dfrac{n\left(\dfrac{b-x}{l_2}\right)-\left(\dfrac{x}{l_1}\right)}{\left(\dfrac{y}{l_1}\right)+\left(\dfrac{ny}{l_2}\right)} $}
\end{equation}
\subsubsection {The Solution of the Differential Equation}

 We use the ``arrow" notation.  ``$P\Rightarrow Q$" means ``the proposition P (logically) implies the proposition Q."

\begin{align*}
\label{}
\eqref{ODE1}  \Rightarrow\left(\dfrac{yy'}{l_1}\right)+\left(\dfrac{nyy'}{l_2}\right)-n\left(\dfrac{b-x}{l_2}\right)-\left(\dfrac{x}{l_1}\right)&=0   \\
  \Rightarrow \dfrac{yy'+x}{l_1}+\dfrac{nyy'-n(b-x)}{l_2}  &=0\\
   \Rightarrow \dfrac{yy'+x}{\sqrt{x^2+y^2}}+n\cdot \dfrac{-(b-x)+yy'}{\sqrt{(b-x)^2+y^2}}  &=0\\
      \Rightarrow\dfrac{1}{2} \dfrac{2yy'+2x}{\sqrt{x^2+y^2}}+n\cdot\dfrac{1}{2} \dfrac{-2(b-x)+2yy'}{\sqrt{(b-x)^2+y^2}}  &=0\\
      \Rightarrow\dfrac{d}{dx}\left\{\sqrt{x^2+y^2}+n\sqrt{(b-x)^2+y^2}\right\}&=0\\
        \Rightarrow \sqrt{x^2+y^2}+n\sqrt{(b-x)^2+y^2}&=c\\
\end{align*}for some (arbitrary) constant $c$.  We have therefore proved:

\begin{thm} The general solution for the differential equation \eqref{ODE1} of the perfect two-dimensional lens, $\mathbf{C}$, with respect to the points $F=(0,0)$ and $F'=(b,0)$, where $b\geqslant 0$, is given by the equation:
\begin{equation}
\label{OVAL}
\fbox{$\dis \sqrt{x^2+y^2}+n\sqrt{(b-x)^2+y^2}=c.$}
\end{equation}
\end{thm}
\hfill$\Box$

As we saw, this is the equation \eqref{car1} of a \textbf{\emph{cartesian oval}} with \textbf{\emph{foci}} at the points $(0,0)$ and $(b,0).$

We have assumed that $b$ is \emph{finite} in this analysis, i.e., that the two foci are a finite distance apart.

Now we consider the limiting cases where one or both foci are ``at infinity."  We will see that we obtain proper or degenerate \emph{\textbf{conic sections}} for these cases.

\subsection{One Fixed Point At Infinity; One Fixed Point Finite}
\subsubsection{The Differential Equation}

We will slightly alter the treatment for the case of two finite foci.  To do so, we begin with the following:

\begin{deff}
A \textbf{\emph{point at infinity}} is specified by means of a line through the origin.  The \emph{\textbf{line joining $P$ to a point at infinity}} is the line through $P$ parallel to the given line.  Points at infinity are not considered to be on $\mathbf{C}$.
\end{deff}

We assume that the fixed point $F$ is at $-\infty$ along the $x$-axis and that the fixed point $F'$ is at the point $(b,0)$ of the $x$-axis, where $b\geqslant 0$.  

Intuitively, this means that \emph{a beam of light from $-\infty$, parallel to the $x$-axis,  is brought to a point focus at $(b,0)$ by a single refracting curve, $f(x,y)=0$, of index $n$.}


The line joining $P$ to the point at infinity is the line parallel to the $x$-axis through $P$.  $\theta_1$ is the angle the horizontal line through $P(x,y)$ makes with the normal while $\theta_2$ is the angle $\overline{PF}$ makes with the normal.  Finally, $l$ be the length of $\overline{PF}$.

Then, the earlier derivation of the ODE is applicable.  We need only observe that $$\theta_1+\phi=\frac{\pi}{2}.$$  So,  substituting our two new formulas for $\sin \theta_1$ and $\sin \theta_2$ into  \textsc{Snell}'s \emph{\textbf{Law}} gives us, instead of \eqref{Snell1},  the new equation:
$$\dfrac{ \dfrac{1}{\sqrt{1+y'^2}}}{{\dfrac{1}{\sqrt{1+y'^2}}\left(\dfrac{b-x}{l}\right)-\dfrac{y'y}{l\sqrt{(1+y'^2)}}}}=n.$$

 After some rearrangement,  we obtain \textbf{\emph{the differential equation of the curve}}:  
\begin{equation}
\label{ODE2}
\fbox{$\dis1-\frac{(b-x)-yy'}{\sqrt{(b-x)^2+y^2}}\cdot n=0.$}
\end{equation}

\subsubsection{The Solution of the Differential Equation}
 The ODE \eqref{ODE2} can be solved by the same computations as we did for the ODE  \eqref{ODE1} which lead us to the 
 
 
 \begin{thm} The general solution for the differential equation \eqref{ODE2} of the perfect two-dimensional lens, $\mathbf{C}$, with the radiant point at $-\infty$ is given by the equation:
\begin{equation}
\label{OVAL1}
\fbox{$\dis x+n\sqrt{(b-x)^2+y^2}=c.$}
\end{equation}where $c$ is an arbitrary constant.
\end{thm}
\hfill$\Box$

We observe that the equation \eqref{OVAL1} has the following interesting interpretation.   The equation  \eqref{OVAL1} says that the ratio of the distance of the point $P$ from the \emph{line} $x=\dfrac{c-bn}{1-n}$ to its distance from the \emph{point} $(b,0)$ is the constant $\pm n$, and thererefore, by the focus-directrix definition, is a \emph{conic section}.
 
 This theorem takes a more elegant form if we assume that the curve $\mathbf{C}$ passes through the origin.  Then, the constant $c=nb$ and, after rationalizing \eqref{OVAL1},  we obtain (\cite{LV}, problem B-10, Chapter 20):
\begin{thm}  The general solution for the differential equation \eqref{ODE2} of the perfect two-dimensional lens, $\mathbf{C}$, is a \emph{\textbf{conic section}} whose \emph{\textbf{focus}} is the point where the light is focused and whose \emph{\textbf{excentricity}} is the \emph{\textbf{reciprocal of the index of refraction}}.   \begin{enumerate}   \item If $n^2\neq 1$, $\mathbf{C}$ given by the equation:
\begin{equation}
\label{CONIC }
\fbox{$\dis \frac{\left(x-\frac{nb}{n+1}\right)^2}{\left(\frac{nb}{n+1}\right)^2}+\frac{y^2}{b^2\left(\frac{n-1}{n+1}\right)}=1$}
\end{equation}Therefore $\mathbf{C}$ is an \emph{\textbf{ellipse}} if $n^2>1$ or an \emph{\textbf{hiperbola}} if $n^2<1$,  either one of which is centered at $\left(\dfrac{nb}{n+1},0\right)$.
  \item If $n=1$, then $\mathbf{C}$ is the \textbf{\emph{segment}} of the $x$-axis given by $0\leqslant x\leqslant b.$
  \item    If $n=-1$, then  $\mathbf{C}$ is the \textbf{\emph{parabola}}
\begin{equation}
\label{PARABOLA }
\fbox{$\dis y^2=4bx$}
\end{equation}

\end{enumerate}\end{thm}
\hfill $\Box$
  
  The reader should compare this result with that of the form of the\textbf{ \emph{perfect reflecting}} \emph{\textbf{mirror}} already cited in \cite{Dru}.  If $n<0$, then we get \textbf{reflection} instead of refraction.
  
    \textsc{Maesumi} \cite{Mae} used \textsc{Fermat}'s Principle to treat this case in a very elegant paper, although his definition of the index of refraction is the reciprocal of our (standard) one.
\subsection{Both Fixed Points are at Infinity}

Keeping the notation of the case of the radiant point at $-\infty$, we assume that the refracted rays form a parallel beam in the direction such that $$\theta_2+\phi=\text{Constant},$$but, this means that$$\sin(\theta_2+\phi)=\frac{b-x}{\sqrt{(b-x)^2+y^2}}=C$$where $C$ is some constant.  But the condition that $\mathbf{C}$ goes through the origin means that $$C=1,$$ and rationalizing the resulting equation we obtain:
\begin{thm}If both fixed points are at infinity, then the perfect lens $\mathbf{C}$ has the equation:\begin{equation}
\label{bothinfinity}
\fbox{$\dis x=0$}
\end{equation}That is, it is the vertical $y$-axis.\end{thm}

\hfill $\Box$

\section{Descartes' Theorem}

\subsection{\textsc{Drucker}'s Characterization of a Surface of Revolution}

In 1992 \cite{Dru} \textsc{Drucker} published a very interesting paper in which he treated the problem of finding all \emph{\textbf{perfect mirrors}}, i.e., mirrors which reflect all rays issuing from one radiant point to one image point.

After showing that the two dimensional curve with the perfect reflecting property is a proper or degenerate conic section, he (implicitly) proved the following characterization of a surface of revolution.  \textsc{Drucker}, himself, did not state it explicitly.
\begin{thm} Let $F$ and $F'$ be two fixed points.  If, for each point $P$ of the smooth connected surface $\mathbf{S}$ the normal $\vec{N}$ at $P$ lies in the subspace spanned by the vectors $\overrightarrow{FP}$ and $\overrightarrow{F'P}$, then $\mathbf{S}$ is a \textbf{surface of revolution} whose axis of revolution is the line through $F$ and $F'$.\end{thm}


\begin{proof}

We offer a \emph{\textbf{new proof}} of \textsc{Drucker}'s theorem.  It is based on an idea in \textsc{Salmon} \cite{Sal} which goes back to \textsc{Monge} \cite{Mon}.

Since, by definition, the normal $\overline{MN}$ is in the subspace spanned by $\overline{FP}$ and $\overline{F'P}$, it is in the plane of $\triangle{FPF'}$.

If $\overline{MN}$ is always parallel to $\overline{FF'}$, then $\mathbf{S}$ is a \textbf{\emph{plane}} which is perpendicular to $\overline{FF'}$.  We exclude this degenerate case for the rest of the argument. (See \eqref{bothinfinity}).

Therefore $\overline{MN}$ is not always parallel to $\overline{FF'}$.  Thus, the infinite line $\overline{MN}$ intersects $\overline{FF'}$ at some point. \emph{This is the characteristic property of the surface} $\mathbf{S}$.

Let $(\alpha,\beta, \gamma)$ be a point on $\overline{FF'}$ and let $(l,m, n)$ be the line's direction numbers where we assume $l\cdot m\cdot n\neq 0$.  The corollaries deal with the case where one or more coefficients are equal to zero.   Then, the equation of the line $\overline{FF'}$ is\begin{equation}
\label{AXIS}
\fbox{$\dis \frac{x-\alpha}{l}=\frac{y-\beta}{m}=\frac{z-\gamma}{n}=t$}
\end{equation}where $t$ is the common value of the three fractions.

Let
\begin{equation}
\label{SURFACE}
F(x,y,z)=0
\end{equation}be the equation of $\mathbf{S}$, where $F$ is a continuously differentiable function of $x$, $y$, and $z$ in some open set $R$, and let $(x_0, y_0,z_0)$ be the point $P$ on $\mathbf{S}$..

Since $\overline{MN}$ is normal to $\mathbf{S}$ at $P$, its equation is:\begin{equation}
\label{NORMAL}
\fbox{$\dis \frac{x-x_0}{F_{x}(x_0, y_0, z_0)}=\frac{y-y_0}{F_{y}(x_0, y_0, z_0)}=\frac{z-z_0}{F_{z}(x_0, y_0, z_0)}=T$}
\end{equation}where $T$ is the common value of the three fractions, and where $F_{x}(x_0, y_0, z_0)\equiv \dfrac{\partial F}{\partial x}$ evaluated in $(x_0,  y_0,z_0)$, and where the other denominators have a similar interpretation.  We assume that all three denominators are different from zero.  The corollaries deal with the cases where the denominators are equal to zero.

Solving equations \eqref{AXIS} and \eqref{NORMAL} for $x$, $y$, and $z$, and then equating the values obtained, we get the following \emph{homogeneous} linear system for the unknowns $t$, $T$, and $1$:\begin{align*}
\label{}
    lt-F_{x}(x_0, y_0, z_0)T+(\alpha-x_0)\cdot 1&=0   \\
    mt-F_{y}(x_0, y_0, z_0)T+(\beta-y_0)\cdot 1&=0  \\
    nt-F_{z}(x_0, y_0, z_0)T+(\gamma-z_0)\cdot 1&=0
\end{align*}The analytical condition that this system have a nontrivial solution, which it \emph{does} by assumption, is that the determinant of their coefficients vanish:\begin{equation}
\label{jacobian}
\begin{vmatrix}
     F_{x}(x_0, y_0, z_0) &F_{y}(x_0, y_0, z_0)&F_{z}(x_0, y_0, z_0)    \\
      l&m&n\\
        x_0-\alpha&y_0-\beta&z_0-\gamma
\end{vmatrix}=0
\end{equation}The determinant on the left-hand side of \eqref{jacobian} is (one half of ) the \textsc{Jacobi}an of the three functions\begin{align}
\label{JACOBIANS}
  \Omega:= F(x,y,z),& &u:=lx+my+nz,& &v:=(x-\alpha)^2+(y-\beta)^2+(z-\gamma)^2,   
      \end{align}
evaluated at the point $P$ of $S$.

But the point $P$ is totally arbitrary, which means the \textsc{Jacobi}an  \eqref{jacobian} vanishes in a full neighborhood of $P$, since $F$ is a continuously differentiable function of  $x$, $y$, and $z$ in some open set $R$.  According to a classical theorem (see \textsc{Buck} \cite{Buc}, \textsc{Goursat} \cite{Gou}, \textsc{Osgood} \cite{Osg}, and \textsc{Taylor} \cite{Tay}, if the \textsc{Jacobi}an of the three functions vanishes identically, then three functions are \emph{\textbf{functionally dependent}}.  

That means that there is a function, $\Omega (u,v)$, of the two variables $u$ and $v$, defined and continuously differentiable in a neighborhood of the point $(u_0,v_0)$, where  $$u_0:=lx_0+my_0+nz_0,\  \ v_0:=(x_0-\alpha)^2+(y_0-\beta)^2+(z_0-\gamma)^2$$ for which the equation\begin{equation}
\label{LENS}
\fbox{$\dis F(x,y,z)=\Omega\left\{lx+my+nz, (x-\alpha)^2+(y-\beta)^2+(z-\gamma)^2\right\}$}
\end{equation}holds identically in a neighborhood of $(x_0,y_0,z_0)$.

Now, the equation\begin{equation}
\label{PLANE}
\fbox{$\dis lx+my+nz=u$}
\end{equation}represents a \emph{plane} which cuts the line $\overline{FF'}$ (represented by \eqref{AXIS}) \emph{perpendicularly}, while the equation\begin{equation}
\label{SPHERE}
\fbox{$\dis (x-\alpha)^2+(y-\beta)^2+(z-\gamma)^2=v$}
\end{equation}represents a \emph{sphere} of radius $\sqrt{v}$ and with center $(\alpha,\beta,\gamma)$ on the line $\overline{FF'}$.

The points $(x,y,z)$ which are on the plane, \eqref{PLANE}, and on the sphere, \eqref{SPHERE},   \emph{simultaneously}, are on their \emph{circle} of intersection and this circle has its center on the line $\overline{FF'}$.

Therefore, the equation \eqref{SURFACE} of $\mathbf{S}$, i.e., $\Omega(u,v)=0$,  represents a surface \emph{generated by a circle of variable radius whose center moves along the line $\overline{FF'}$ and whose plane is perpendicular to that line.}

Thus, every planar transverse section of $\mathbf{S}$, perpendicular to $\overline{FF'}$, consists of one or more circles whose centers are on the line $\overline{FF'}$.  

That is, $\mathbf{S}$ is a \emph{surface of revolution with axis $\overline{FF'}$.}

This completes the proof of \textsc{Drucker}'s theorem.
\end{proof}


\begin{cor}If the $z$-axis is the axis of revolution, we may take the origin as the point $(\alpha,\beta,\gamma)$, and the equation \eqref{SURFACE} becomes \begin{equation}
\label{LENS1}
\fbox{$\dis F(x,y,z)=\Omega\left\{z, x^2+y^2+z^2\right\}$}
\end{equation}

\end{cor}

\hfill $\Box$

There are similar simplifications in \eqref{LENS} if we take the other coordinate axes as the axis of revolution.

\begin{cor}If  $F_{x}(x_0, y_0, z_0)\equiv 0$, then the normal is everywhere perpendicular to the $x-axis$  and the equation \eqref{SURFACE} becomes the \textbf{cylinder} of revolution: \begin{equation}
\label{LENS2}
\fbox{$\dis F(x,y,z)=\Omega\left\{my+nz, (y-\beta)^2+(z-\gamma)^2\right\}$}
\end{equation}

\end{cor}

\hfill $\Box$

There are similar simplifications if the other components of the normal are zero.

\begin{rmk}

\emph{In order to apply the theorem on \emph{\textbf{functional dependence}} which we used in the above proof, we have to make sure that we comply with all the hypotheses.  The only one, which we did not explicitly state in the body of the proof is, using the notations of \eqref{PLANE} and \eqref{SPHERE}, is that} at least one of the three jacobians\begin{align}
\label{JACOBIANS}
    \frac{\partial(u,v)}{\partial(x,y)},& &\frac{\partial(u,v)}{\partial(y,z)},& &\frac{\partial(u,v)}{\partial(x,z)},   
      \end{align}is different from zero at $(x_0,y_0,z_0)$.
       
\end{rmk}We claim that  \emph{even more is true} in our case.  We will prove that \emph{at least \textbf{two} of the jacobians \eqref{JACOBIANS} are different from zero.} 

Suppose, to the contrary, that at least two of them are equal to zero, say
\begin{align}
\label{JACOBIANS1}
    \frac{\partial(u,v)}{\partial(x,y)}&=0 &\frac{\partial(u,v)}{\partial(y,z)}&=0   
      \end{align} This leads to\begin{align}
\label{JACOBIANS1}
    \frac{l}{x-\alpha}&=\frac{m}{y-\beta}, &\frac{m}{y-\beta}&=\frac{n}{z-\gamma},   
      \end{align}respectively.  By \eqref{JACOBIANS1}
      \begin{align}
\label{JACOBIANS2}
    \frac{l}{x-\alpha}=\frac{m}{y-\beta} &=\frac{n}{z-\gamma}   
      \end{align}which is the equation of the \emph{axis} $\overline{FF'}$.  But, this means that $\mathbf{S}$ is just the \emph{straight line axis}, which is excluded by the hypothesis that $\mathbf{S}$ is a smooth surface.  Therefore, at least two of the jacobians \eqref{JACOBIANS} are different from zero and the theorem on functional dependence is applicable.
      
\begin{rmk}
\emph{The proof shows that the characteristic property of a surface $\mathbf{S}$ of revolution is that} the normal to any point of $\mathbf{S}$ intersects the axis of revolution.      \end{rmk}
\subsection{Proof of Descartes' Theorem}

We adapt \textsc{Drucker}'s definition 

\begin{deff} Let $\mathbf{S}$ be a smooth connected surface and let $F$ and $F'$ be points not in $\mathbf{S}$.  We say that $\mathbf{S}$ is a \textbf{\emph{perfect lens}} relative to $F$ and $F'$ if, for each point $P$ in $\mathbf{S}$:\begin{enumerate}
  \item  the normal $\vec{N}$ at $P$ lies in the subspace spanned by the vectors $\overrightarrow{FP}$ and $\overrightarrow{F'P}$, and
  \item the sines of the angles which  $\overrightarrow{FP}$ and $\overrightarrow{F'P}$ form with that normal are in \emph{\textbf{constant ratio}} for every point $P$ in $\mathbf{S}$.
  \end{enumerate}
\end{deff}

 By condition $2$ of the definition, the cross-section of $\mathbf{S}$ sliced out by the $xy$-plane is a plane curve $\mathbf{C}$ which is a perfect two-dimensional lens relative to $F$ and $F'$.
 
 That means that $\mathbf{C}$ is either (part of) a \emph{\textbf{cartesian oval}}, or (part of) a \emph{\textbf{conic section}}, or a degenerate case of either one.
 
 Therefore, by condition $1$ and \textsc{Drucker}'s Theorem, a three dimensional perfect lens $\mathbf{S}$ is (part of) a surface of \emph{\textbf{revolution}} with axis $\overline{FF'}$ obtained by rotating a two-dimensional perfect lens $\mathbf{S}$ around it.

This completes the proof of  \textsc{Descartes}' Theorem.

\hfill $\Box$

\subsubsection*{Acknowledgment}
 Support from the Vicerrector\'{\i}a de Investigaci\'on of the 
University of Costa Rica is acknowledged.

\end{document}